\documentclass[12pt]{amsart}
\usepackage{amssymb}
\usepackage[margin=1in]{geometry}
\usepackage[all]{xy}

\long\def\forget#1\forgotten{}
\newcommand{\nc}{\newcommand}
\nc{\my}[1]{\textsf{(#1)}\marginpar{**}}
\nc{\PR}{\op{PR}}
\nc{\seq}[2]{\langle\, #1 : #2\,\rangle}
\nc{\bbN}{\mathbb{N}}
\nc{\weight}{\op{w}}
\nc{\density}{\op{d}}
\nc{\localdensity}{\op{ld}}
\nc{\boundedness}{\op{b}}
\nc{\Nh}{\widehat{\cN}}
\nc{\la}{\langle}
\nc{\ra}{\rangle}
\nc{\myqn}[1]{{\sf #1}\marginpar{???}}
\nc{\cl}[1]{\overline{#1}}
\nc{\inv}{^{-1}}
\nc{\kw}{\op{kw}}
\nc{\w}{\omega}
\nc{\Fin}[1]{\op{Fin}(#1)} 
\nc{\finw}[1]{{\Fin{#1}^\bbN}}
\nc{\kfinw}{\finw{\kappa}}
\nc{\fb}{\mathfrak{b}}
\nc{\fd}{\mathfrak{d}}
\nc{\fc}{\mathfrak{c}}
\nc{\ff}{\mathfrak{f}}
\nc{\Bgp}{\mathbb{Z}^\bbN}
\nc{\mult}{*}
\nc{\Cite}[1]{\textbf{[#1]}}
\nc{\Pa}[9]{\bibitem{#1} {#2}, \emph{#3}, {#4} \textbf{#5} ({#6}), {#7}--{#8}.}
\nc{\Bc}[9]{\bibitem{#1} {#2}, \emph{#3}, in: \textbf{#4} (#5), #6 #7, #8--#9.}
\nc{\alephes}{{\aleph_0}}
\nc{\oo}{\infty}
\nc{\beq}{\begin{eqnarray*}}
\nc{\eeq}{\end{eqnarray*}}
\nc{\op}[1]{\operatorname{#1}}
\nc{\PK}{\op{PK}}
\nc{\K}{\op{K}}
\nc{\rmPK}{\mathrm{PK}}
\nc{\cov}{\op{cov}}
\nc{\non}{\op{non}}
\nc{\cof}{\op{cof}}
\nc{\Par}{\op{Par}}
\nc{\vphi}{\varphi}
\nc{\Bdd}{\op{Bdd}}
\nc{\card}[1]{\left|#1\right|}
\nc{\AP}{\mathrm{A\!P}}
\nc{\RH}{\op{R\!H}}
\nc{\FP}{\op{FP}}
\nc{\FS}{\op{FS}}
\nc{\comp}{^\mathsf{c}}
\nc{\Impl}{\Rightarrow}
\nc{\bi}{\begin{itemize}}
\nc{\itm}{\item}
\nc{\ei}{\end{itemize}}
\nc{\be}{\begin{enumerate}}
\nc{\ee}{\end{enumerate}}
\nc{\Tau}{\mathrm{T}}
\nc{\nin}{\notin}
\nc{\cO}{\mathcal{O}}
\nc{\Un}{\bigcup}
\nc{\sub}{\subseteq}
\nc{\sps}{\supseteq}
\nc{\bbQ}{\mathbb{Q}}
\nc{\NN}{{\bbN^\bbN}}
\nc{\Z}{\mathbb{Z}}
\nc{\bbC}{\mathbb{C}}
\nc{\bbT}{\mathbb{T}}
\nc{\bbB}{\mathbb{B}}
\nc{\R}{\mathbb{R}}
\nc{\set}[2]{\{\, #1 : #2\,\}}
\nc{\Icov}[1]{{#1}_{\cI}}
\nc{\scrA}{\mathscr{A}}
\nc{\scrB}{\mathscr{B}}
\nc{\scrC}{\mathscr{C}}
\nc{\sm}{\setminus}
\nc{\FinSeqs}[1]{{#1}^*}
\nc{\as}{\subseteq^*}
\nc{\roth}{{[\bbN]^{\!\/\infty}}}
\nc{\setroth}[1]{{[#1]^{\!\/\infty}}}
\nc{\grp}[1]{\gimel(#1)}
\nc{\grpa}{\grp{\scrA}}
\nc{\grpb}{\grp{\scrB}}
\nc{\sone}{\mathsf{S}_{1}}
\nc{\sfin}{\mathsf{S}_\mathrm{fin}}
\nc{\ufin}{\mathsf{U}_\mathrm{fin}}
\nc{\Split}{\mathsf{Split}}
\nc{\gone}{\mathsf{G}_{1}}
\nc{\goab}{\gone(\mathcal{A},\scrB)}
\nc{\gfin}{\mathsf{G}_\mathrm{fin}}
\nc{\cA}{\mathcal{A}}
\nc{\cM}{\mathcal{M}}
\nc{\cH}{\mathcal{H}}
\nc{\cI}{\mathcal{I}}
\nc{\OI}{\cO_\cI}
\nc{\cB}{\mathcal{B}}
\nc{\cX}{\mathcal{X}}
\nc{\cU}{\mathcal{U}}
\nc{\cV}{\mathcal{V}}
\nc{\cW}{\mathcal{W}}
\nc{\cK}{\mathcal{K}}
\nc{\cF}{\mathcal{F}}
\nc{\cG}{\mathcal{G}}
\nc{\cS}{\mathcal{S}}
\nc{\cR}{\mathcal{R}}
\nc{\cP}{\mathcal{P}}
\nc{\cN}{\mathcal{N}}
\nc{\x}{\times}
\nc{\st}{\op{star}}

\newtheorem{thm}{Theorem}[section]
\nc{\bthm}{\begin{thm}} \nc{\ethm}{\end{thm}}
\newtheorem{prop}[thm]{Proposition}
\nc{\bprp}{\begin{prop}} \nc{\eprp}{\end{prop}}
\newtheorem{fact}[thm]{Fact}
\nc{\bfct}{\begin{fact}} \nc{\efct}{\end{fact}}
\newtheorem{prob}[thm]{Problem}
\nc{\bprb}{\begin{prob}} \nc{\eprb}{\end{prob}}
\newtheorem{lem}[thm]{Lemma}
\nc{\blem}{\begin{lem}} \nc{\elem}{\end{lem}}
\newtheorem{claim}[thm]{Claim}
\nc{\bclm}{\begin{claim}} \nc{\eclm}{\end{claim}}
\newtheorem{cor}[thm]{Corollary}
\nc{\bcor}{\begin{cor}} \nc{\ecor}{\end{cor}}
\newtheorem{conj}[thm]{Conjecture}
\nc{\bcnj}{\begin{conj}} \nc{\ecnj}{\end{conj}}
\theoremstyle{definition}
\newtheorem{defn}[thm]{Definition}
\nc{\bdfn}{\begin{defn}} \nc{\edfn}{\end{defn}}
\newtheorem{spec}[thm]{Specializing}
\nc{\bspc}{\begin{spec}} \nc{\espc}{\end{spec}}
\theoremstyle{remark}
\newtheorem{rem}[thm]{Remark}
\nc{\brem}{\begin{rem}} \nc{\erem}{\end{rem}}
\newtheorem{cnv}[thm]{Convention}
\nc{\bcnv}{\begin{cnv}} \nc{\ecnv}{\end{cnv}}
\newtheorem{exam}[thm]{Example}
\nc{\bexm}{\begin{exam}} \nc{\eexm}{\end{exam}}
\nc{\bpf}{\begin{proof}} \nc{\epf}{\end{proof}}

\nc{\ed}{

\end{document}
}

\title[Star covering properties]{Combinatorial aspects of selective star covering properties in $\Psi$-spaces}

\author{Boaz Tsaban}
\address{Department of Mathematics, Bar-Ilan University, Ramat Gan 5290002, Israel,
and Faculty of Mathematics and Computer Science, Weizmann Institute of Science, Rehovot 76100\-01, Israel}
\email{tsaban@math.biu.ac.il}
\urladdr{math.biu.ac.il/~tsaban}

\subjclass[2010]{
54D20, 
54A35, 
03E17. 
}

\keywords{Isbell--Mr\'owka space, $\Psi$-space, star-Menger, star-Hurewicz.}

\begin{document}

\begin{abstract}
Which Isbell--Mr\'owka spaces ($\Psi$-spaces)
satisfy the star version of Menger's and Hurewicz's covering properties?
Following Bonanzinga and Matveev, this question is considered here from a combinatorial point of view.
An example of a $\Psi$-space that is (strongly) star-Menger but not star-Hurewicz is obtained.
The PCF-theory function $\kappa\mapsto\cof([\kappa]^\alephes)$ is a key tool.
Using the method of forcing, a complete answer to a question of Bonanzinga and Matveev is provided.

The results also apply to the mentioned covering properties in the realm of
Pixley--Roy spaces, to the extent of spaces with these properties, and to
the character of free abelian topological groups over hemicompact $k$ spaces.
\end{abstract}

\maketitle

\section{Introduction}

The Isbell--Mr\'owka $\Psi$-spaces~\cite{Isbell, Mrowka54}
are classic examples in the realm of topological covering properties.
A family $\cA\sub P(\bbN)$ is \emph{almost disjoint} if every element of $\cA$ is infinite, and the sets
$A\cap B$ are finite for all distinct elements $A,B\in \cA$.
For an almost disjoint family $\cA$, let $\Psi(\cA):=\cA\cup\bbN$. A topology on $\Psi(\cA)$ is defined as follows.
The natural numbers are isolated, and for each element $A\in\cA$ and each finite set $F\sub\bbN$,
the set $\{A\}\cup (A\sm F)$ is a basic open neighborhood of $A$. Spaces constructed in this manner
are called \emph{$\Psi$-spaces}.

For a set $X$, a subset $A$ of $X$ and a family $\cU$ of subsets of $X$, let
$\st(A,\cU):=\Un\set{U\in\cU}{A\cap U\ne \emptyset}$.
A topological space $X$ is \emph{star-Lindel\"of}~\cite{star} if
every open cover $\cU$ of $X$ has a countable subset $\cV$ such that $X=\st(\Un\cV,\cU)$.
It is \emph{strongly star-Lindel\"of}~\cite{star} if,
for each open cover $\cU$ of $X$,
there is a countable set $C\sub X$ such that $X=\st(C,\cU)$.
It is easy to see that uncountable $\Psi$-spaces are not Lindel\"of. Being separable, though,
all $\Psi$-spaces are strongly star-Lindel\"of.

\emph{Menger's property} is the following selective version of Lindel\"of's property:
For every sequence $\cU_1,\cU_2,\dots$ of open covers of $X$,
there are finite sets $\cF_1\sub\cU_1, \cF_2\sub\cU_2, \dots$ such that the family 
$\{\Un\cF_1,\Un\cF_2,\dots\}$
covers $X$.

A topological space $X$ is
\emph{star-Menger}
(respectively, \emph{strongly star-Menger})
\cite{Koc99} if for every sequence $\cU_1,\cU_2,\dots$ of open covers of $X$,
there are finite sets $\cF_1\sub\cU_1, \cF_2\sub\cU_2, \dots$
(respectively, $F_1,F_2,\cdots\sub X$)
such that the family $\{\st(\Un\cF_1,\cU_1),\st(\Un\cF_2,\cU_2),\dots\}$
(respectively, $\{\st(F_1,\cU_1),\st(F_2,\cU_2),\dots\}$) covers $X$.

A topological space $X$ is a \emph{Hurewicz} (respectively: \emph{star-Hurewicz};
\emph{strongly star-Hurewicz}) \emph{space}~\cite{BCK04} if,
in the corresponding definitions in the previous paragraph, we request that every point of $X$ is in
the set $\Un\cF_n$ (respectively: $\st(\Un\cF_n,\cU_n)$; $\st(F_n,\cU_n)$) for all but finitely many $n$.

The implications among the mentioned covering properties are as follows. 
$$
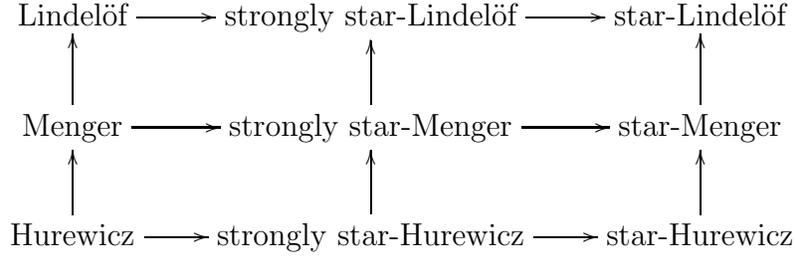
\xymatrix{
\mbox{Lindel\"of}\ar[r] & \mbox{strongly star-Lindel\"of}\ar[r] & \mbox{star-Lindel\"of}\\
\mbox{Menger}\ar[r]\ar[u] & \mbox{strongly star-Menger}\ar[r]\ar[u] & \mbox{star-Menger}\ar[u]\\
\mbox{Hurewicz}\ar[r]\ar[u] & \mbox{strongly star-Hurewicz}\ar[r]\ar[u] & \mbox{star-Hurewicz}\ar[u]
}$$
A survery of these properties and their connections to other notions is available in~\cite{KocPsiSurv}.

Background on the combinatorial cardinals of the continuum used in this paper, including the unbounding number
$\fb$ and the dominating number $\fd$, is available in~\cite{vD, BlassHBK}. Whether a $\Psi$-space
is strongly star-Menger---or strongly star-Hurewicz---depends only on the cardinality of the space.

\bthm[Bonanzinga--Matveev~\cite{MilenaMisha}]\label{thm:ssMH}
Let $\cA\sub P(\bbN)$ be an almost disjoint family.
\be
\item The space $\Psi(\cA)$ is strongly star-Menger if and only if $\card{\cA}<\fd$.
\item The space $\Psi(\cA)$ is strongly star-Hurewicz if and only if $\card{\cA}<\fb$.
\ee
\ethm

The question of when a $\Psi$-space $\Psi(\cA)$ is star-Menger---or star-Hurewicz---is more elusive.
Combinatorial characterizations in terms of the family $\cA$
are provided in Section~\ref{sec:cc}, but some of the most basic problems remain, in general, open.
Some of these problems are reviewed in Section~\ref{sec:probs}.

Let $P$ be a partially ordered set. A subset $C$ of $P$ is \emph{cofinal} if
for each element $a\in P$ there is an element $c\in C$ such that $a\le c$.
The \emph{cofinality} of $P$, denoted $\cof(P)$, is the minimal cardinality
of a cofinal subset of $P$. The number $\cof(P)$ may, in general, be a singular cardinal number.
For a set $X$, let $\Fin{X}$ be the family of all finite subsets of $X$.
In this paper, families of sets are always partially ordered by the relation $\sub$.
The set $\Fin{X}^\bbN$ of all functions $f\colon\bbN\to\Fin{X}$ is partially ordered
coordinate-wise: $f\le g$ if $f(n)\sub g(n)$ for all $n$.
The cardinal $\cof(\Fin{X}^\bbN)$ depends only on $\card{X}$.
For an infinite cardinal $\kappa$, the cardinal $\cof(\finw{\kappa})$
will later be expressed in simpler terms. In particular, it is known that
the cardinality $\fc$ of the continuum satisfies $\cof(\finw{\fc})=\fc$.

\bthm[Bonanzinga--Matveev~\cite{MilenaMisha}]\label{thm:NSM}
Let $\cA\sub P(\bbN)$ be an almost disjoint family of cardinality $\kappa$.
If $\cof(\finw{\kappa})=\kappa$, then the space $\Psi(\cA)$ is not star-Menger.
\ethm

A simple proof of Theorem~\ref{thm:NSM} is provided in Section~\ref{sec:cc}.
Section~\ref{sec:cc} also includes a similar theorem for star-Hurewicz $\Psi$-spaces (Theorem~\ref{thm:NSH}).
Theorems~\ref{thm:ssMH}(1) and~\ref{thm:NSH} are used in Example~\ref{exm:SSMnotSH} to obtain
a consistent example of a (strongly) star-Menger $\Psi$-space that is not star-Hurewicz.

The existence of a star-Menger $\Psi$-space that is not star-Hurewicz violates the 
Continuum Hypothesis, and thus cannot be constructed in ZFC alone.
Indeed, $\Psi$-spaces have cardinality at most $\fc$. Since $\cof(\finw{\fc})=\fc$, 
every star-Menger $\Psi$-space has cardinality smaller than $\fc$. 
By Theorem~\ref{thm:ssMH}(2), we have the following corollary.

\bcor\label{cor:Ari}
If $\fb=\fc$, then every star-Menger $\Psi$-space is (strongly) star-Hurewicz.\qed
\ecor

\brem
If we do not insist on $\Psi$-spaces then there is, provably in ZFC,
a very nice (strongly) star-Menger space that is not star-Hurewicz:
For paracompact spaces, each of the mentioned covering properties coincides with its star- and strongly star-
versions. Chaber and Pol proved that there are Menger subsets of the Cantor space
that are not Hurewicz (cf.~\cite{MHP}).
\erem

The question whether $\cof(\finw{\kappa})=\kappa$ for a cardinal number $\kappa$
appears in a number of additional, related and seemingly unrelated, topological contexts.
The following theorem follows from Sakai's Theorem~2.1 in~\cite{Sakai14}, since
being closed discrete is a hereditary property.

\bthm[Sakai]\label{cor:ff}
Let $D$ be a closed discrete subspace of a regular strongly star-Menger space.
Then the cardinality of $D$ is smaller than the minimal fixed point of the function $\kappa\mapsto\cof(\finw{\kappa})$.
\ethm

Let $X$ be a topological space. The \emph{Pixley--Roy space} $\PR(X)$ is the space of all nonempty finite
subsets of $X$, with the topology determined by the basic open sets
$$[F,U] := \set{H\in\PR(X)}{F\sub H\sub U},$$
$F\in\PR(X)$ and $U$ open in $X$.

\bthm[Sakai~\cite{Sakai14}]\label{thm:s}
Let $X$ be an infinite regular topological space of cardinality $\kappa$. If $\cof(\finw{\kappa})=\kappa$,
then the space $\PR(X)$ is not star-Menger.
\ethm

The cardinals $\cof(\Fin{\kappa}^\bbN)$ also show up in a study of the character of topological groups.

\bthm[\cite{PvK}]
Let $X$ be a nondiscrete hemicompact $k$ space.
Let $\kappa$ be the supremum of the weights of compact subsets of $X$.
Then the character of the free abelian topological group $A(X)$
is $\cof(\Fin{\kappa}^\bbN)$.
\ethm

A similar result is proved in~\cite{PvK} for general abelian non-locally compact hemicompact $k$ groups.
A number of estimations of $\cof(\Fin{\kappa}^\bbN)$ for infinite cardinals $\kappa$ are provided there.
The key to these is the following reduction. For an infinite cardinal number $\kappa$, let $[\kappa]^\alephes$
be the family of all countably infinite subsets of $\kappa$.

\bprp[\cite{PvK}]\label{morph}
Let $\kappa$ be an infinite cardinal number. Then
$\cof(\Fin{\kappa}^\bbN)$ is the maximum of the cardinals $\fd$ and $\cof([\kappa]^\alephes)$.
\eprp

Thus, the Bonanzinga--Matveev Theorem~\ref{thm:NSM} can be reformulated as follows.
(Recall that the space $\Psi(\cA)$ is strongly star-Menger if $\card{\cA}<\fd$.)

\bthm\label{thm:often}
Let $\cA\sub P(\bbN)$ be an almost disjoint family of cardinality $\kappa\ge\fd$.
If $\cof([\kappa]^\alephes)=\kappa$, then the space $\Psi(\cA)$ is not star-Menger. 
\ethm

The estimation of the cardinal $\cof([\kappa]^\alephes)$ in terms of the cardinal $\kappa$
is a central goal in Shelah's \emph{PCF theory}, the theory of possible cofinalities.
In contrast to cardinal exponentiation, the function $\kappa\mapsto\cof([\kappa]^\alephes)$
is tame. For example, if there are no large cardinals in the Dodd--Jensen core model,
then $\cof([\kappa]^\alephes)$ is simply $\kappa$ if $\kappa$ has uncountable cofinality,
and $\kappa^+$ (the successor of $\kappa$) otherwise~\cite{Gitik91}.
Moreover, without any special hypotheses, the cardinal
$\cof([\kappa]^\alephes)$ can be estimated, and in many cases computed exactly.
Some examples follow (for proofs and references, see~\cite[Section 8]{PvK}).

For uncountable cardinals $\kappa$ of countable cofinality, a variation of K\"onig's Lemma implies that $\cof([\kappa]^\alephes)>\kappa$.
Throughout, \emph{Shelah's Strong Hypothesis (SSH)} is the assertion that
$\cof([\kappa]^\alephes)=\kappa^+$ for all uncountable cardinals $\kappa$ of countable cofinality.
Clearly, the Generalized Continuum Hypothesis implies SSH, but the latter axiom is much 
weaker, being a consequence of the absence of large cardinals.

\bthm[Folklore]\label{nice}
The following cardinals are fixed points of the function $\kappa\mapsto\cof\allowbreak([\kappa]^\alephes)$:
\be
\item The cardinals $\kappa$ with $\kappa^\alephes=\kappa$.
\item $\aleph_n$, for natural numbers $n\ge 1$.
\item The cardinals $\aleph_\kappa$, for $\kappa$ a singular cardinal of uncountable cofinality
that is smaller than the first fixed point of the $\aleph$ function.
\item Assuming SSH, all cardinals of uncountable cofinality.
\ee
Moreover, successors of fixed points of this function are also fixed points.
\ethm

For example, for $n=1,2,\dots$, the cardinal $\aleph_{\aleph_{\w_n}}$ and its successors are all fixed
points of the function $\kappa\mapsto\cof([\kappa]^\alephes)$.

\bcor
Let $\cA\sub P(\bbN)$ be an almost disjoint family of cardinality at least $\fd$.
\be
\itm
For each cardinal $\kappa$ smaller than the first fixed point of
the $\aleph$ function, with $\alephes<\cof(\kappa)<\kappa$,
if $\card{\cA}=\aleph_\alpha$ for some ordinal $\alpha$ with $\kappa\le\alpha<\kappa+\omega$,
then the space $\Psi(\cA)$ is not star-Menger.
\itm Assume SSH. If the cardinal $\card{\cA}$ has uncountable cofinality, then
the space $\Psi(\cA)$ is not star-Menger.\qed
\ee
\ecor

The cardinality of $\Psi$-spaces is at most $\fc$.
Knowing that $\cof(\Fin{\kappa}^\bbN)=\fd\cdot\kappa$ for the cardinals $\aleph_n$ (for $n\in\bbN$)
and for the cardinal $\fc$, the following problem is natural.

\bprb[Bonanzinga--Matveev~\cite{MilenaMisha}]\label{MM}
Is $\cof(\kfinw)=\fd\cdot\kappa$ for each infinite cardinal $\kappa\le\fc$?
In particular, is $\cof(\kfinw)=\fd$ for each infinite cardinal $\kappa\le\fd$?
\eprb

This problem is solved in Section~\ref{sec:pcfcon}.

\section{Combinatorial characterizations and a consequence}\label{sec:cc}

The following theorem provides a combinatorial characterization of star-Menger $\Psi$-spaces.

\bthm\label{thm:char}
Let $\cA\sub P(\bbN)$ be an almost disjoint family.
The following assertions are equivalent:
\be
\item The Isbell--Mr\'owka space $\Psi(\cA)$ is star-Menger.
\item For each function $A\mapsto f_A$ from $\cA$ to $\NN$,
there are finite sets $\cF_1,\cF_2,\dots\sub\cA$ such that, for each $A\in\cA$,
there is $n$ with $(A\sm f_A(n))\cap \Un_{B\in\cF_n}(B\sm f_B(n))\neq\emptyset$.
\ee
\ethm
\bpf
$(2)\Impl (1)$:
Since the subspace $\bbN$ of $\Psi(\cA)$ is countable, it suffices in the definition of the star-Menger property
to cover $\cA$.
Let $\cU_n$, for $n\in\bbN$, be open covers of $\Psi(\cA)$.
By moving to a finer open cover, we may assume that
for each $A\in\cA$ and each $n$, there is a natural number $f_A(n)$
such that $\{A\}\cup(A\sm f_A(n))\in\cU_n$.

Let $\cF_1,\cF_2,\dots\sub\cA$ be finite sets as in (2).
For each $n$, the set
$$\set{\{B\}\cup(B\sm f_B(n))}{B\in\cF_n}$$
is a finite subset of $\cU_n$.
Let $A\in\cA$. Pick $n$ as in (2).
Then
$$A\in\{A\}\cup (A\sm f_A(n))\sub \st\bigl(\,\mbox{$\Un$}_{B\in\cF_n}(\{B\}\cup(B\sm f_B(n))),\cU_n\,\bigr).$$

$(1)\Impl (2)$:
For each $n$, let
$$\cU_n:=\set{ \{A\}\cup (A\sm f_A(n)) }{ A\in\cA } \cup \set{\{m\} }{ m\in\bbN}.$$
Since the space $\Psi(\cA)$ is star-Menger, there are finite sets $\cF_1\sub\cU_1,\cF_2\sub\cU_2,\dots$
such that $\Psi(\cA)=\Un_n\st(\Un\cF_n,\cU_n)$.
For each $n$ and each $\{m\}\in\cF_n$, pick if possible an element $B\in\cA$ such that
$m\in B\sm f_B(n)$, and substitute $\{B\}\cup (B\sm f_B(n))$ for $\{m\}$ in $\cF_n$.
If there is no such $B$, just remove $\{m\}$ from $\cF_n$ (in this case, $\st(\{m\},\cU_n)=\{m\}$).
Then $\cA\sub\Un_n\st(\Un\cF_n,\cU_n)$. The assertion in (2) then follows from the definitions.
\epf

We obtain the following simple proof of Theorem~\ref{thm:NSM}.
The main simplification over the proof in~\cite{MilenaMisha}
is that we avoid the necessity to use two types of cofinal sets simultaneously.

\bpf[Proof of Theorem~\ref{thm:NSM}]
We establish the negation of the characterization in Theorem~\ref{thm:char}.

Enumerate $\cA:=\set{A_\alpha}{\alpha<\kappa}$, and let $\set{F_\alpha}{\alpha<\kappa}$ be a cofinal
subset of $\Fin{\kappa}^\bbN$.
We may assume that $\alpha\notin F_\alpha(n)$ for all $n$. Indeed,
the family $\set{F_\alpha'}{\alpha<\kappa}$, defined by $F_\alpha'(n):=F_\alpha(n)\sm\{\alpha\}$ 
for all $n$, is cofinal in $\Fin{\kappa}^\bbN$:
Let $F\in\kfinw$, and set $I:=\set{\alpha<\kappa}{F\le F_\alpha}$.
For each ordinal $\beta<\kappa$, there is $\alpha<\kappa$ such that
$F(n)\cup\{\beta\}\sub F_\alpha(n)$ for all $n$.
Thus, $\Un_{\alpha\in I}\Un_n F_{\alpha}(n)=\kappa$,
and therefore the set $I$ is uncountable. 
Pick an ordinal $\alpha\in I\sm \Un_nF(n)$.
Then $F(n)\sub F_\alpha(n)\sm\{\alpha\}$ for all $n$. 

For each $\alpha<\kappa$ and each $n$, let
$$f_{\alpha}(n):=1+\max\Un_{\beta\in F_\alpha(n)} A_\alpha\cap A_\beta.$$
Let $\cF_1,\cF_2,\dots\sub\cA$ be finite sets. For each $n$, let
$H_n:=\set{\alpha<\kappa}{A_\alpha\in\cF_n}$.
Take $\alpha$ such that $H_n\sub F_\alpha(n)$ for all $n$.
Then, for each $n$, we have that
$\max (A_\alpha\cap \Un_{\beta\in H_n}A_\beta)
<f_\alpha(n)$,
and thus
$$
(A_\alpha\sm f_\alpha(n))\cap \Un_{\beta\in H_n}(A_\beta\sm f_\beta(n))\sub
(A_\alpha\sm f_\alpha(n))\cap \Un_{\beta\in H_n}A_\beta=\emptyset.\qedhere$$
\epf

The following theorem provides a combinatorial characterization of star-Hurewicz $\Psi$-spaces.
Its proof, which is similar to that of Theorem~\ref{thm:char}, is omitted.

\bthm\label{thm:charH}
Let $\cA\sub P(\bbN)$ be an almost disjoint family.
The following assertions are equivalent:
\be
\item The Isbell--Mr\'owka space $\Psi(\cA)$ is star-Hurewicz.
\item For each function $A\mapsto f_A$ from $\cA$ to $\NN$,
there are finite sets $\cF_1,\cF_2,\dots\sub\cA$ such that, for each $A\in\cA$,
$(A\sm f_A(n))\cap \Un_{B\in\cF_n}(B\sm f_B(n))\neq\emptyset$ for all but finitely many $n$.\qed
\ee
\ethm

\bprp\label{morph2}
Let $\kappa$ be an infinite cardinal.
The following cardinal numbers are equal:
\be
\item The minimal cardinality of a family $\cF\sub \Fin{\kappa}^\bbN$ such that for each $g\in\Fin{\kappa}^\bbN$
there is $f\in\cF$ with $g(n)\sub f(n)$ for infinitely many $n$.
\item The maximum of the cardinals $\fb$ and $\cof([\kappa]^\alephes)$.
\ee
\eprp
\bpf
$(2)\le (1)$: Let $\cF$ be as in (1).

For each $f\in\cF$, define a function $f'\in\NN$ by
$$f'(n):=1+\max (f(n)\cap\bbN).$$
For each function $g\in\NN$, there is $f\in\cF$ such that $\{1,\dots,g(n)\}\sub f(n)$, and thus
$g(n) \le f'(n)$, for infinitely many $n$. Thus, the family $\set{f'}{f\in\cF}$ is unbounded. This shows that
$\fb\le\card{\cF}$.

For each set $A\in[\kappa]^\alephes$, pick a function $g\in\kfinw$ such that $g(n)\sub g(n+1)$ for all $n$,
and $A\sub \Un_ng(n)$. Pick $f\in\cF$ such that $g(n)\sub f(n)$ for infinitely many $n$.
Then, since $g(n)\sub g(n+1)$ for all $n$, $\Un_n g(n)\sub \Un_n f(n)$.
Thus, the family $\set{\Un_n f(n)}{f\in\cF}$ is cofinal in $[\kappa]^\alephes$. It follows that
$\cof([\kappa]^\alephes)\le\card{\cF}$.

$(1)\le (2)$: Let $\cG$ be an unbounded family in $\NN$, and $\cH$ be a cofinal family in $[\kappa]^\alephes$.
For each set $A\in\cH$, fix a function $f_A\in\kfinw$ such that $f_A(n)\sub f_A(n+1)$ for all $n$,
and $A\sub \Un_nf_A(n)$.

Let $h\in\kfinw$. Pick $A\in\cH$ with $\Un_nh(n)\sub A$. Pick $g\in\cG$
such that
$$\min\set{m}{h(n)\sub f_A(m)}\le g(n)$$
for infinitely many $n$. Then $h(n)\sub f_A(g(n))$ for infinitely many $n$.
Take $\cF:=\set{ f_A\circ g}{g\in\cG, A\in\cH}$.
Then $\card{\cF}\le\fb\cdot \cof([\kappa]^\alephes)$.
\epf

We obtain the following analogue of Theorem~\ref{thm:often}. 
(Recall that $\Psi$-spaces of cardinality smaller than $\fb$ are strongly star-Hurewicz.)

\bthm\label{thm:NSH}
Let $\cA\sub P(\bbN)$ be an almost disjoint family of cardinality $\kappa\ge\fb$.
If $\cof([\kappa]^\alephes)=\kappa$, then the space $\Psi(\cA)$ is not star-Hurewicz.
\ethm
\bpf
The proof is almost identical to 
that of Theorem~\ref{thm:NSM}, using Proposition~\ref{morph2} and
Theorem~\ref{thm:charH}.
The necessary changes are as follows.
Here, we let $\set{F_\alpha}{\alpha<\kappa}\sub\kfinw$ be a family as in Proposition~\ref{morph2}(1).
For the last step of the proof, we take $\alpha$ such that $H_n\sub F_\alpha(n)$ for infinitely many $n$,
and restrict attention to these $n$.
\epf

\bexm\label{exm:SSMnotSH}
Assume that $\fb=\aleph_1<\fd$. Then there is a strongly star-Menger $\Psi$-space that is not star-Hurewicz.
\eexm
\bpf
Since there are almost disjoint sets of cardinality continuum, there are ones of any smaller cardinality, too.
Let $\cA\sub P(\bbN)$ be an almost disjoint family of cardinality $\aleph_1$.
By Theorem~\ref{thm:ssMH}, the space $\Psi(\cA)$ is strongly star-Menger.
By Theorem~\ref{nice}(2) and Theorem~\ref{thm:NSH}, this space is not star-Hurewicz.
\epf

\bcor[SSH]
The following assertions are equivalent:
\be
\item There is a strongly star-Menger $\Psi$-space that is not star-Hurewicz.
\item $\fb<\fd$.
\ee
\ecor
\bpf
$(1)\Impl (2)$: Let $\Psi(\cA)$ exemplify (1).
By Theorem~\ref{thm:ssMH}, $\card{\cA}<\fd$. If $\fb=\fd$ then,
by the same theorem, the space $\Psi(\cA)$ is (strongly) star-Hurewicz; a contradiction.

$(2)\Impl (1)$: Take a $\Psi$-space of cardinality $\fb$.
By Theorem~\ref{thm:ssMH}, the space $\Psi(\cA)$ is strongly star-Menger.
By Theorem~\ref{nice}(4), since $\fb$ is a regular cardinal, $\cof([\fb]^\alephes)=\fb$.
Apply Theorem~\ref{thm:NSH}.
\epf

\section{A solution of the Bonanzinga--Matveev Problem}\label{sec:pcfcon}

Problem~\ref{MM} asks whether $\cof(\kfinw)=\fd\cdot\kappa$ for each infinite cardinal $\kappa\le\fc$, 
and, in particular, whether $\cof(\kfinw)=\fd$ for each infinite cardinal $\kappa\le\fd$.

Clearly, the Continuum Hypothesis implies a positive answer to Problem~\ref{MM},
and Problem~\ref{MM} actually asks whether the assertions are provable without special
set theoretic hypotheses. 
We first point out a negative answer to the first part of this problem.

\bprp\label{prp:simple}
Let $\aleph_\alpha:=\fd$. If $\aleph_{\alpha+\omega}<\fc$, 
then there is a cardinal $\kappa<\fc$ such that $\cof(\kfinw)>\fd\cdot\kappa$.
\eprp
\bpf
Take $\kappa:=\aleph_{\alpha+\omega}$. Since $\fd\le\kappa\le \cof([\kappa]^\alephes)$, 
we have by Theorem~\ref{nice} that $\cof(\kfinw)=\cof([\kappa]^\alephes)$.
By K\"onig's Lemma, we have that $\cof([\kappa]^\alephes)>\kappa=\fd\cdot\kappa$.
\epf

We use some facts from the theory of forcing. A general introduction is available in Kunen's book~\cite{KunenST}, whose notation
we follow.
Some more details that are relevant for us here are available in Bartoszy\'nski and Judah's book~\cite{barju},
and in Blass's chapter~\cite{BlassHBK}.

Fix a successor ordinal $\beta>\omega$.
Adding $\aleph_\beta$ random reals to a model of the Continuum Hypothesis,
we obtain a model of $\fd=\aleph_1$ and $\fc=\aleph_\beta$. Such a model satisfies the condition
in Proposition~\ref{prp:simple}.

SSH implies a positive answer to the second part of the Bonanzinga--Matveev problem,
and a conditional solution to its first part.

\bthm[SSH]\label{MMSSH}
\mbox{}
\be
\itm For each infinite cardinal $\kappa\le\fd$, we have that $\cof(\kfinw)=\fd$.
\itm $\cof(\kfinw)=\fd\cdot\kappa$ for all infinite cardinals $\kappa\le\fc$ if, and only if,
there is $n\ge 0$ such that $\fc=\fd^{+n}$, the $n$-th successor of $\fd$.
\ee
\ethm
\bpf
We use Theorem~\ref{nice}.

(1) If $\cof(\kappa)>\alephes$, then $\cof(\kfinw)=\fd\cdot\kappa=\fd$.
Otherwise, as $\cof(\fd)\ge\fb>\alephes$, we have that $\kappa<\fd$, and
$\cof(\kfinw)=\fd\cdot\kappa^+=\fd$.

(2) If there is such $n$, then each $\kappa$ with $\fd\le\kappa\le\fc$ has uncountable cofinality,
and by SSH we have that $\cof(\kfinw)=\fd\cdot\kappa$.
Otherwise, Proposition~\ref{prp:simple} applies.
\epf

Thus, the answer to the first part of Problem~\ref{MM} is ``No'', and
the answer to its second part is ``Yes'' if there are no (inner) models of set theory with large cardinals.
To complete the picture, it remains to show that the answer is ``No'' (to both parts)
when large cardinal hypotheses are available. For the following theorem, it suffices for example to
assume the consistency of supercompact cardinals, or of so-called \emph{strong cardinals}. More
precise large cardinal hypotheses are available in~\cite{GitikMagidor92}.

\bthm[Gitik--Magidor~\cite{GitikMagidor92}]\label{MagThm}
It is consistent (relative to the consistency of ZFC with an appropriate large cardinal hypothesis) that
$2^{\aleph_n}=\aleph_{n+1}$ for all $n$, and $2^{\aleph_\omega}=\aleph_{\omega+\gamma+1}$,
for any prescribed $\gamma<\omega_1$.
\ethm

This theorem is related to our questions as follows. As $\aleph_\w$ is a limit cardinal of cofinality $\alephes$,
$2^{\aleph_\omega}=(2^{<\aleph_\w})^\alephes$. If $2^{\aleph_n}=\aleph_{n+1}$ for all $n$,
then $2^{<\aleph_\w}=\aleph_\w$, and thus $2^{\aleph_\omega}=(\aleph_\w)^\alephes=2^\alephes\cdot\cof([\aleph_\w]^\alephes)=
\cof([\aleph_\w]^\alephes)$.\footnote{For the second equality, count the countable subsets of $\aleph_\w$
by taking a cofinal family in $[\aleph_\w]^\alephes$ and, for each set in this family, take all of its subsets.}

\emph{Hechler's forcing} $\mathbb{H}$ is a natural forcing notion adding a dominating real, i.e.,
$d\in\NN$ such that for each $f\in\NN\cap V$, where $V$ is the ground model, $f\le^* d$.
$\mathbb{H}=\set{(n,f)}{n\in\bbN, f\in\NN}$, and $(n,f)\le(m,g)$ if $n\ge m$, $f\ge g$, and $f(k)=g(k)$ for all $k<m$.
If $G$ is $\mathbb{H}$-generic over $V$, then by a density argument, $d=\Un_{(n,f)\in G}f|_{\{1,\dots,n\}}\in\NN$
is as required. $\mathbb{H}$ is ccc, and thus so is the finite support iteration
$P=\seq{P_\alpha,\dot Q_\alpha}{\alpha<\lambda}$,
where for each $\alpha$, $P_\alpha$ forces that $\dot Q_\alpha$ is Hechler's forcing.

\bthm\label{ans}
It is consistent (relative to the consistency of ZFC with appropriate large cardinal hypotheses) that
$$\aleph_\w<\fb=\fd=\aleph_{\w+1}<\cof(\finw{\aleph_\w})=\cof([\aleph_\w]^\alephes)=\aleph_{\w+\gamma+1}=\fc,$$
for each prescribed $\gamma$ with $1\le\gamma<\aleph_1$.
\ethm
\bpf
Use Theorem~\ref{MagThm} to produce a model of set theory, $V$, satisfying $\fc=\aleph_1$ and $\cof([\aleph_\w]^\alephes)\allowbreak=\aleph_{\w+\gamma+1}$.

Let $P:=\seq{P_\alpha,\dot Q_\alpha}{\alpha<\aleph_{\w+1}}$ be the
finite support iteration, where for each $\alpha$, $P_\alpha$ forces that $\dot Q_\alpha$ is Hechler's forcing.
Let $G$ be $P$-generic over $V$, and for each $\alpha<\aleph_{\w+1}$, let $G_\alpha:=G\cap P_\alpha$ be the induced
$P_\alpha$-generic filter over $V$. For each $\alpha$, let $d_\alpha$ be the dominating real added by $Q_\alpha$ in
stage $\alpha+1$, so that for each $f\in V[G_\alpha]\cap\NN$, $f\le^* d_\alpha$.

As $P$ is ccc, $\cof([\aleph_\w]^\alephes)$ remains $\aleph_{\w+\gamma+1}$ in $V[G]$.
As $\aleph_{\w+1}$ has uncountable cofinality, we have that $\NN\cap V[G]=\Un_{\alpha<\aleph_{\w+1}}\NN\cap V[G_\alpha]$
\cite[Lemma 1.5.7]{barju}.
It follows that $\set{d_\alpha}{\alpha<\aleph_{\w+1}}$ is dominating in $V[G]$. Moreover, it follows that
for each $B\sub\NN\cap V[G]$ with $\card{B}<\aleph_{\w+1}$, there is $\alpha<\aleph_{\w+1}$ such that $B\sub\NN\cap V[G_\alpha]$,
and thus $B$ is $\le^*$-bounded (by $d_\alpha$). Thus, in $V[G]$, $\fb=\fd=\aleph_{\w+1}$.

As the Continuum Hypothesis holds in $V$, $\card{P}=\aleph_{\w+1}$, and
as $P$ is ccc, the value of $\fc$ in $V[G]$ is at most (by counting nice names~\cite[Lemma 5.13 in Chapter VII]{KunenST})
$\card{P}^\alephes=\aleph_{\w+1}^{\alephes}$, evaluated in $V$.
In $V$, $\aleph_{\w+1}^{\alephes}\le (2^{\aleph_\w})^\alephes=2^{\aleph_\w}=\aleph_{\w+\gamma+1}$.
Thus, in $V[G]$, $\fc\le\aleph_{\w+\gamma+1}$.
On the other hand, in $V[G]$, as $\aleph_\w<\fd\le\fc$, $\aleph_{\w+\gamma+1}=\cof([\aleph_\w]^\alephes)\le\aleph_\w^\alephes\le\fc^\alephes=\fc$.
\epf

\brem\label{rem:Rinot}
For finite $\gamma$, which are sufficient for our purposes, a simplified proof of the Gitik--Magidor Theorem~\ref{MagThm}
is available in Gitik's chapter~\cite{GitikHBK}.
Following our proof, Assaf Rinot pointed out to us that
starting with a supercompact cardinal (a stronger assumption than that in~\cite{GitikHBK}),
one may argue as follows:
Start with a model of GCH with $\kappa$ supercompact.
Use Silver forcing to make $2^\kappa=\kappa^{++}$~\cite[Theorem 21.4]{JechST}.
Since $\kappa$ remains measurable, we can use Prikry forcing to make $\cof(\kappa)=\alephes$,
without adding bounded subsets~\cite[Theorem 21.10]{JechST}.
Then GCH holds up to $\kappa$, and $\cof([\kappa]^\alephes)=\kappa^\alephes=2^\kappa=\kappa^{++}$.
Then, continue as in the proof of Theorem~\ref{ans}.
\erem

\section{Comments and open problems}\label{sec:probs}

Remarkably, the following problem remains open.

\bprb[Bonanzinga--Matveev~\cite{MilenaMisha}]\label{prb:main}
Is there, consistently, a star-Menger $\Psi$-space of cardinality $\ge\fd$?
\eprb

Since $\Psi$-spaces of cardinality smaller than $\fd$ are strongly star-Menger, the problem
asks whether there could be star-Menger $\Psi$-spaces that are not in fact strongly star-Menger.
More importantly, the problem asks whether there may be, consistently, \emph{nontrivial} 
star-Menger $\Psi$-spaces,
that is, ones whose being star-Menger does not follow from their cardinality being smaller than $\fd$.
By Theorem~\ref{thm:often},
the cardinality of a nontrivial star-Menger $\Psi$-space cannot be any of the cardinals listed in Theorem~\ref{nice}.
Thus, $\fc>\aleph_\omega$ in every model witnessing a positive solution of Problem~\ref{prb:main}.
It may be worth considering forcing extensions where 
$\fd=\aleph_1$, $\kappa=\aleph_\omega$, and $\fc=\aleph_{\omega+1}$.
Similarly, we have the following problem (to which similar comments apply).

\bprb
Is there, consistently, a star-Hurewicz $\Psi$-space of cardinality $\ge\fb$?
\eprb

A topological space $X$ is \emph{star-Rothberger}~\cite{Koc99} if 
for every sequence $\cU_1,\cU_2,\dots$ of open covers of $X$,
there are elements $U_1\in\cU_1, U_2\in\cU_2,\dots$ such that $X=\Un_n\st(U_n,\cU_n)$.
Arguments similar to ones in Section~\ref{sec:cc} establish the following theorem.

\bthm\label{thm:charR}
Let $\cA\sub P(\bbN)$ be an almost disjoint family.
The following assertions are equivalent:
\be
\item The Isbell--Mr\'owka space $\Psi(\cA)$ is star-Rothberger.
\item For each function $A\mapsto f_A$ from $\cA$ to $\NN$,
there are elements $A_1,A_2,\dots\in\cA$ such that, for each $A\in\cA$,
there is $n$ with $(A\sm f_A(n))\cap (A_n\sm f_{A_n}(n))\neq\emptyset$.\qed
\ee
\ethm

The cardinal $\cov(\cM)$ is the minimal cardinality of a subset of $\NN$ that cannot be guessed
by a single function (that is, no function is equal infinitely often to each member of the set).
It is open whether there is an analogue of Theorems~\ref{thm:often} and~\ref{thm:NSH}
for star-Rothberger $\Psi$-spaces. $\Psi$-spaces of cardinality smaller than $\cov(\cM)$ are
star-Rothberger, and there is $\Psi$-space of cardinality $\cov(\cM)$ that is not star-Rothberger~\cite{MilenaMisha}.

\bprb
Is there, consistently, an almost disjoint family $\cA\sub P(\bbN)$ of cardinality $\kappa\ge\cov(\cM)$
such that $\cof([\kappa]^\alephes)=\kappa$ and the space $\Psi(\cA)$ is star-Rothberger?
\eprb

It is not clear that the cardinals in Theorems~\ref{thm:NSM} and~\ref{thm:NSH} are not
mere artifact of the proofs. Indeed, the proofs exploit the freedom provided by
Theorems~\ref{thm:char} and~\ref{thm:charH}.
In particular, we have the following problems.

\bprb
What is the minimal cardinal $\kappa$ such that no $\Psi$-space of cardinality $\kappa$ is star-Menger?
What is the corresponding cardinal for star-Hurewicz and star-Rothberger $\Psi$-spaces?
\eprb

In light of Section~\ref{sec:cc}, it may be possible to prove, using the methods of~\cite{Sakai14}, 
the following variations of Theorems~\ref{cor:ff} and~\ref{thm:s} 

\bcnj
\mbox{}
\be
\item Let $D$ be a closed discrete subspace of a regular strongly star-Hurewicz space.
Then the cardinality of $D$ is smaller than the minimal fixed point of the function
$\kappa\mapsto\cof([\kappa]^\alephes)$ in the interval $[\fb,\fc]$.

\item Let $X$ be a regular topological space of cardinality $\kappa$. If $\cof([\kappa]^\alephes)=\kappa\ge\fb$,
then the space $\PR(X)$ is not star-Hurewicz.
\ee
\ecnj

Motivated by Theorem~\ref{cor:ff}, Sakai proposes the following problem.

\bprb[Sakai]
Consider the minimal cardinal number greater than all cardinalities of
closed discrete subspaces of regular strongly star-Menger spaces. Is this cardinal
equal to the minimal fixed point of the function $\kappa\mapsto\cof([\kappa]^\alephes)$ in the interval
$[\fd,\fc]$?
\eprb

\subsection*{Acknowledgments}
I thank Moti Gitik for bringing Theorem~\ref{MagThm} to my attention
and for his useful suggestions,
Assaf Rinot for his comment in Remark~\ref{rem:Rinot},
and Masami Sakai and Shir Sivroni for their useful comments.
I owe special thanks to Ari Meir Brodsky, whose comments helped improving the 
presentation of this paper considerably, and to the referee for a detailed and useful 
report. 

A part of the research reported here was conducted during a Sabbatical leave
at the Faculty of Mathematics and Computer Science, Weizmann Institute of Science. 
I thank Gideon Schechtman and the Faculty of Mathematics and Computer Science for their hospitality.
\ed